\documentclass[11pt]{article}

\usepackage{enumerate}
\usepackage{amsmath}
\usepackage{amssymb,latexsym}
\usepackage{amsthm}

\usepackage{anysize}
\usepackage{graphicx}
\usepackage{wrapfig}

\date{}

\newcommand{\cN}{{\mathcal N}}
\newcommand{\cC}{{\mathcal C}}
\newcommand{\cM}{{\mathcal M}}

\newcommand{\cL}{{\mathcal L}}

\newcommand{\GF}{\hbox{{\rm GF}}}
\newcommand{\st}{{\cal S}_t} 
\newcommand{\so}{{\cal S}} 
 
\newcommand{\bo}{{\cal B}}

\renewcommand{\proof}{\noindent{\bf Proof.}\ }
\newcommand{\Qed}{\hfill $\Box$ \medskip}

\newcommand{\PG}{\mathrm{PG}}
\newcommand{\AG}{\mathrm{AG}}

\newtheorem{theorem}{Theorem}[section]

\newtheorem{lemma}[theorem]{Lemma}
\newtheorem{corollary}[theorem]{Corollary}
\newtheorem{definition}[theorem]{Definition}
\newtheorem{proposition}[theorem]{Proposition}

\newtheorem{example}[theorem]{Example}

\begin{document}

\title{Semiarcs with long secants}
\author{Bence Csajb\'{o}k\thanks{Author was supported by the Hungarian National Foundation for Scientific Research, Grant No. K 81310.}}

\maketitle

\begin{abstract}
In a projective plane $\Pi_q$ of order $q$, a non-empty point set $\st$ is a $t$-semiarc if the number of tangent lines to $\st$ at each of its points is $t$.
If $\st$ is a $t$-semiarc in $\Pi_q$, $t<q$, then each line intersects $\st$ in at most $q+1-t$ points.
Dover proved that semiovals (semiarcs with $t=1$) containing $q$ collinear points exist in $\Pi_q$ only if $q<3$.
We show that if $t>1$, then $t$-semiarcs with $q+1-t$ collinear points exist only if $t\geq \sqrt{q-1}$.
In $\PG(2,q)$ we prove the lower bound $t\geq(q-1)/2$, with equality only if $\st$ is a blocking set of R\'edei type of size $3(q+1)/2$.

We call the symmetric difference of two lines, with $t$ further points removed from each line, a $V_t$-configuration.
We give conditions ensuring a $t$-semiarc to contain a $V_t$-configuration and give the complete characterization of such $t$-semiarcs in $\PG(2,q)$.

%%%%%%%%%%%%%%%%%%%%%%%%%%%%%%%%%%%%%%%%%%%%%%%%%%%%%%%%%%%%%%%%%%%%%%%%%%%%%%%%%%%%%%%%%%%%%%%%%%%%%%%%%%%%%%%%%%%%%%%%%%%%%%%%%%%%%%%%%%%%%%%%%%%%%

%A $t$-semiarc of a finite projective plane is a pointset $\st$ with the property that the number of tangent lines to $\st$ at each of its points is $t$.
%When $t=1$, then $\st$ is called a semioval. If $\st$ is a $t$-semiarc in a projective plane of order $q$ with $q>t$, then each line intersects $\st$ in at most $q+1-t$ points. We generalize a result of Dover about semiovals with a long secant and prove that if $t>1$, then $t$-semiarcs with $q+1-t$ collinear points exist only if $t\geq \sqrt{q-1}$. Using a conjecture of Metsch, that was proved by Weiner and Sz\H{o}nyi, we prove the lowerbound $t\geq(q-1)/2$ in $\PG(2,q)$,
%which is sharp as the projective triangle shows.
%
%We call the symmetric difference of two lines, with $t$ further points removed from each line a $V_t$-configuration.
%We give conditions ensuring a $t$-semiarc to contain a $V_t$-configuration and give the complete characterization of such $t$-semiarcs in $\PG(2,q)$.
%In the proof we use the classification of perspective pointsets in $\PG(2,q)$.
%This is a result due to Korchm\'{a}ros and Mazzoca and it is related to Dickson's classification of the subgroups of the affine group on the line $\AG(1,q)$.
\end{abstract}

\section{Introduction}
\label{sec:Intro}
\indent

Semiarcs are natural generalizations of arcs. Let $\Pi _q$ be a projective plane of order $q$.
A non-empty point set $\st\subset \Pi _q$ is called a {\it t-semiarc} 
if for every point $P\in \st$ there exist exactly $t$ lines $\ell _1,\ell _2,\ldots, \ell _t$ such that 
$\st \cap \ell_i = \{P\}$ for $i=1,2,\ldots ,t$.
These lines are called the {\it tangents} to $\st$ at $P$. If a line $\ell $ meets $\st$ in $k>1$ points, then $\ell$ is called
a {\it $k$-secant} of $\st$.
The classical examples of semiarcs are the semiovals (semiarcs with $t=1$) and point sets of type $(0, 1, n)$
(i.e. point sets meeting each line in either 0, or 1, or $n$ points, in this case $t=q + 1 - (s-1)/(n-1)$, where $s$ denotes the size of the point set).
Arcs, unitals and subplanes are semiarcs of the latter type. For more examples, see \cite{2semi}, \cite{smallsemi} and \cite{kgy2}.

Because of the huge diversity of the geometry of semiarcs, their complete classification is hopeless. 
In \cite{dov} Dover investigated semiovals with a $q$-secant and semiovals with more than one $(q-1)$-secant.
The aim of this paper is to generalize these results and characterize $t$-semiarcs with long secants.

Many of the known $t$-semiarcs contain the symmetric difference of two lines, with $t$ further points removed from each line.
We will call this set of $2(q-t)$ points a {\it $V_t$-configuration}.
Recently in $\cite{smallsemi}$ it was proved that in $\PG(2,q)$ small semiarcs with a long secant necessarily contain a $V_t$-configuration
or can be obtained from a blocking set of R\'{e}dei type.
Here we give another condition ensuring a $t$-semiarc to contain a $V_t$-configuration and we give the complete characterization of such $t$-semiarcs in $\PG(2,q)$. To do this we use the classification of perspective point sets in $\PG(2,q)$.
This is a result due to Korchm\'{a}ros and Mazzoca \cite{nuclei} and it is related to Dickson's classification of the subgroups of the affine group on the line $\AG(1,q)$.

Using a result of Weiner and Sz\H{o}nyi, that was conjectured by Metsch, we prove that $t$-semiarcs in $\PG(2,q)$ with $q+1-t$ collinear points
exist if and only if $t\geq (q-1)/2$. The case of equality is strongly related to blocking sets of R\'{e}dei type, we also
discuss these connections.

If $t=q+1,q$ or $q-1$, then $\st$ is single point, a subset of a line or three non-collinear points respectively.
To avoid trivial cases, we may assume for the rest of this paper that $t<q-1$.
%If $t=q-2$, then it follows from \cite{csk}, Proposition 3.1, that 
%$\st$ is one of the following three configurations:  
%four points in general position, the six vertices of a complete 
%quadrilateral, or a Fano subplane. 

\section{Semiarcs with one long secant}
\label{sec:1.1}
\indent

If $\st$ is a $t$-semiarc in $\Pi_q$, $t<q$, then each line intersects $\st$ in at most $q+1-t$ points.
In this section we study $t$-semiarcs containing $q+1-t$ collinear points.
The following lemma gives an upper bound for the size of such $t$-semiarcs.

\begin{lemma}
\label{hosszu}
If $\st$ is a $t$-semiarc in $\Pi_q$ and $\ell$ is a $(q+1-t)$-secant of $\st$, then $|\st \setminus \ell|\leq q$.
\end{lemma}
\proof
Let $U=\st\setminus \ell$ and let $D=\ell \setminus \st$.
Through each point of $U$ there pass exactly $t$ tangents to $\st$ and each of them intersects $\ell$ in $D$.
This implies $t|U|\leq q|D|$.
Since $|D|=t$, we have $|U|\leq q$.
\Qed

In \cite{dov} Dover proved that semiovals with a $q$-secant exist in $\Pi_q$ if and only if $q\leq 3$.
Our first theorem generalizes this result and shows that if $\st$ has a $(q+1-t)$-secant, then $t$ cannot be arbitrary.
For related ideas of the proof, see the survey paper by Blokhuis et. al. \cite{stability}, Theorem 3.2.
%, and Theorem 1.1 of \cite{joco} from Pelik\'{a}n.

\begin{theorem}
\label{i0}
If $\st$ is a $t$-semiarc in $\Pi_q$ with a $(q+1-t)$-secant, then $t=1$ and $q\leq 3$ or $t\geq \sqrt{q-1}$.
\end{theorem}
\proof
Let $\ell$ be a line that satisfies $|\st \cap \ell|=q+1-t$ and let $U=\st \setminus \ell$.
The size of $U$ has to be at least $q-t$, otherwise the points of $\ell \cap \st$ would have more than $t$ tangents.
This and Lemma \ref{hosszu} together yield:
\begin{equation}
q-t \leq |U| \leq q.
\end{equation}
Let $q-t+k$ be the size of $U$, where $0\leq k \leq t$.
Let $\delta$ be the number of lines that do not meet $U$ and denote by $L_1, L_2, \ldots, L_{q^2+q+1-\delta}$ the lines that meet $U$.
For these lines let $e_i=|L_i \cap U|$. The standard double counting argument gives:
\begin{equation}
\sum_{i=1}^{q^2+q+1-\delta} e_i = (q-t+k)(q+1),
\end{equation}
\begin{equation}
\label{kettes}
\sum_{i=1}^{q^2+q+1-\delta} e_i(e_i-1)=(q-t+k)(q-t+k-1).
\end{equation}
If a line $\ell'$ intersects $U$ in more than one point, then $Q:=\ell' \cap \ell$ is in $\st$, otherwise the points of $\ell' \cap U$
would have at most $t-1$ tangents. The point $Q\in \st$ has at least $q-1-(q-t+k-|\ell'\cap U|)=t-1-k+|\ell'\cap U|$ tangents, hence $|\ell'\cap U| \leq k+1$.
This implies $e_i \leq k+1$, for $i=1,2,\ldots,q^2+q+1-\delta$, thus the following holds:
\begin{equation}
\label{harmas}
\sum_{i=1}^{q^2+q+1-\delta} e_i(e_i-1)\leq (k+1)\sum_{i=1}^{q^2+q+1-\delta}(e_i-1)=(k+1)((q-t+k)(q+1)-(q^2+q+1-\delta)).
\end{equation}
The line $\ell$ does not meet $U$ and the other lines that do not meet $U$ fall into two classes:
there are $(q+1-t)t$ of them passing through $\ell \cap \st$ (the tangents to $\st$ through the points of $\ell\cap \st$) and
there are $tq-(q-t+k)t$ of them passing through $\ell \setminus \st$ (the lines intersecting $\ell \setminus \st$ minus the tangents to $\st$ through the points of $U$). This implies $\delta=t(q+1-k)+1$, hence we can write (\ref{harmas}) as:
\begin{equation}
(q-t+k)(q-t+k-1)\leq(k+1)((q-t+k)(q+1)-(q^2+q)+t(q+1-k)).
\end{equation}
Rearranging this inequality we obtain:
$$q^2-q(2t+1-k+k^2)+k^2t-kt-2k+t^2+t\leq 0.$$
The discriminant of the left-hand side polynomial is $k^4-2k^3+3k^2+6k+1$.
If $k=0,1,2$, then we get $q\leq t+1,t+2,t+4$ respectively.
Otherwise we have $k^4-2k^3+3k^2+6k+1 < (k^2-k+3)^2$, which yields $q \leq t+k^2-k+1$.
The maximum value of $k$ is $t$, therefor $q \leq t^2+1$ follows for $k\geq 3$.
If $t=1$, then $k\leq 1$, hence $q\leq t+2=3$.
If $t=2$, then $k\leq 2$, hence $q\leq t+4=6$. Since there is no projective plane of order 6, in this case we get $q\leq 5$.
If $t\geq 3$ and $k<3$, then $q \leq t+4 < t^2+1$ and this completes the proof.
\Qed

Before we go further we need some definitions about blocking sets.
A {\it blocking set} of a projective plane is a point set $\bo$ that intersects every line in the plane.
A blocking set is {\it minimal} if it does not contain a smaller blocking set and it is {\it non-trivial} if it does not contain a line.
If $\bo$ is a non-trivial blocking set, then we have $|\ell\cap \bo|\leq |\bo|-q$ for every line $\ell$.
If there is a line $\ell$ such that $|\ell\cap \bo|=|\bo|-q$, then $\bo$ is a blocking set of {\it R\'{e}dei type} and the line $\ell$
is a {\it R\'{e}dei line} of $\bo$.

In $\PG(2,q)$ we can improve the bound in Theorem \ref{i0}.
To do this we use the following result, that was proved with the resultant method by Weiner and Sz\H{o}nyi in \cite{SzW,WThesis} and was conjectured by Metsch.

\begin{theorem}[\cite{SzW,WThesis}]
\label{lemma0}
Let $U$ be a point set in $\PG(2,q)$, $P$ a point not from $U$ and assume that there pass exactly $r$ lines through $P$ meeting $U$.
Then the total number of lines meeting $U$ is at most $1+rq+(|U|-r)(q+1-r)$.
\end{theorem}

\begin{theorem}
\label{ii1}
Let $\st$ be a $t$-semiarc in $\PG(2,q)$. If $\st$ has a $(q+1-t)$-secant, then $t\geq(q-1)/2$.
In the case of equality $\st$ is a blocking set of R\'{e}dei type and its $(q+1-t)$-secants are R\'{e}dei lines.
\end{theorem}
\proof
Let $\ell$ be a $(q+1-t)$-secant of $\st$ and let $U=\st \setminus \ell$.
From Lemma \ref{hosszu}, we have:
\begin{equation}\label{eq1}
|U|\leq q.
\end{equation}
The following statements are easy to check:
\begin{itemize}
	\item the lines intersecting $U$ in more than one point intersect $\ell$ in $\ell \cap \st$,
	\item through each point of $\ell \cap \st$ there pass exactly $r=q-t$ lines meeting $U$,
	\item the total number of lines meeting $U$ is $\delta=|U|t+(q+1-t)(q-t)$.
\end{itemize}
Applying Theorem \ref{lemma0} for the point set $U$ and for a point $P\in \ell \cap \st$, we obtain:
\begin{equation}
\delta=|U|t+(q+1-t)(q-t) \leq 1 + (q-t)q+(|U|-q+t)(t+1).
\end{equation}
After rearranging, we get:
\begin{equation}\label{eq2}
2q-2t-1 \leq |U|.
\end{equation}
Equations (\ref{eq1}) and (\ref{eq2}) together imply $t\geq (q-1)/2$.
If $t = (q-1)/2$, then $|U|=q$ and there are $\delta=(3q^2+2q+3)/4$ lines meeting $U$ and $(q+1-t)t=(q^2+2q-3)/4$ lines
meeting $\ell \cap \st$ but $U$. Together with the line $\ell$ we get the total number of lines in $\PG(2,q)$,
thus $\st$ is a blocking set of R\'{e}dei type and $\ell$ is a R\'{e}dei line of $\st$.
\Qed

With respect to the other direction we cite the following result by Blokhuis:

\begin{theorem}[\cite{blokp}]
\label{blok}
If $\bo$ is a minimal non-trivial blocking set in $\PG(2,p)$, $p>2$ prime, then $|\bo|\geq 3(p+1)/2$.
In the case of equality there pass exactly $(p-1)/2$ tangent lines through each point of $\bo$.
\end{theorem}

\begin{example}[\cite{JWP}, Lemma 13.6]
\label{pt}
Denote by $C$ the set of non-zero squares in $\GF(q)$, $q$ odd, and let
$\st=\{(c,0,1),(0,-c,1),(c,1,0) : c \in C\}\cup\{(1,0,0),(0,1,0),(0,0,1)\}$.
This point set is called projective triangle and it is a $t$-semiarc with three $(q+1-t)$-secants, where $t=(q-1)/2$.
This example shows the sharpness of Theorems \ref{ii1} and \ref{blok}.
\end{example}

%If $t=(q-1)/2$ and $\st$ is a $t$-semiarc with a $(q+1-t)$-secant denoted by $\ell$, then Theorem \ref{ii1} implies that $|\st \setminus \ell|=q$.
%Consider $\ell$ as $\ell_{\infty}$, that is the line at infinity. The points of $\ell_{\infty}$ are also called directions.
%A direction $P$ is said to be determined by an affine pointset $\so$ if and only if there is a line through $P$ intersecting $S$ in at least two points.
%Lov\'{a}sz and Schrijver proved that if $\so$ is a set of $q$ affine points in $\PG(2,q)$, $q$ prime, and $\so$ determines $(q+3)/2$ directions, then $\so$ is affinely equivalent to the graph of $x^{(q+1)/2}$, see \cite{LS}. G\'{a}cs, Lov\'{a}sz and Sz\H{o}nyi proved that the same follows if $q$ is a square of a prime, see \cite{GLSZ}. The $q$ points of $\st \setminus \ell$ determine exactly $q+1-t=(q+3)/2$ directions, namely the points of $\st \cap \ell$ hence we have the following result:
In $\PG(2,q)$, $q$ prime, Lov\'{a}sz and Schrijver proved that blocking sets of R\'{e}dei type of size $3(q+1)/2$ are projectively equivalent to the projective
triangle, see \cite{LS}. G\'{a}cs, Lov\'{a}sz and Sz\H{o}nyi proved the same if $q$ is a square of a prime, see \cite{GLSZ}.
These results and Theorem $\ref{ii1}$ together yield the following:

\begin{corollary}[\cite{GLSZ,LS}]
Let $\st$ be a $t$-semiarc in $\PG(2,q)$ with a $(q+1-t)$-secant.
If $t=(q-1)/2$ and $q=p$ or $q=p^2$, $p$ prime, then $\st$ is projectively equivalent to the projective triangle.
\end{corollary}

\section{Semiarcs with two long secants}
\label{sec:1.2}
\indent

Throughout the paper, if $A$ and $B$ are two point sets in $\Pi_q$, then $A\triangle B$ denotes their symmetric difference, that is $(A\setminus B)\cup(B\setminus A)$.

\begin{definition}
A $V_t$-configuration is the symmetric difference of two lines, with $t$ further points removed from both lines.
Semiarcs containing a $V_t$-configuration fall into two types.
Let $\st$ be a $t$-semiarc and suppose that there are two lines, $\ell_1$ and $\ell_2$, such that $(\ell_1 \triangle \ell_2)\cap \st$ is
a $V_t$-configuration, then:
\begin{itemize}
	\item $\st$ is of $V_t^\circ$ type if $\ell_1 \cap \ell_2 \notin \st$,
	\item $\st$ is of $V_t^\bullet$ type if $\ell_1 \cap \ell_2 \in \st$.
\end{itemize}
\end{definition}

For semiovals, Dover proved the following characterization:

\begin{theorem}[\cite{dov}, Lemma 4.1, Theorem 4.2]
\label{dover}
Let $\so_1$ be a semioval in $\Pi_q$. If $\so_1$ is of $V_t^\circ$ type, then it is contained in a vertexless triangle.
If $q>5$ and $\so_1$ has at least two $(q-1)$-secants, then $\so_1$ is of $V_t^\circ$ type.
\end{theorem}

As the above result suggests, the characterization of $t$-semiarcs with two $(q-t)$-secants works nicely only for semiarcs of $V_t^\circ$ type.
In Proposition \ref{dovv} we generalize the last statement of the above result, but the characterization
of $V_t^\circ$ type semiarcs seems to be hopeless in general. In Proposition \ref{le2} we consider the case when $t=2$, but
for larger values of $t$ we deal only with the Desarguesian case, see Section \ref{sec:1.3}.

\begin{lemma}
\label{j1}
Let $\st$ be a $t$-semiarc in $\Pi_q$, $t<q$, and suppose that there exist two lines, $\ell_1$ and $\ell_2$, with their common point in $\st$ such that $|\ell_1 \setminus (\st\cup \ell_2)|=n$ and $|\ell_2 \setminus (\st\cup \ell_1)|=m$.
Then $q \leq t+1+nm/t$ and $|\st \setminus (\ell_1 \cup \ell_2)|=q-1-t$ in the case of equality.
\end{lemma}
\proof
Since $\st$ is not contained in a line, we have $n,m\geq t$. If one of $n$ or $m$ is equal to $q$, then
$q < q+t+1 \leq t+1+nm/t$ and the assertion follows. Thus we can assume that $\ell_1$ and $\ell_2$ are not tangents to $\st$.
Let $X=\st \setminus (\ell_1 \cup \ell_2)$. Through the point $\ell_1 \cap \ell_2 $ there pass exactly $t$ tangents to $\st$, hence $q-1-t \leq |X|$.
Through the points of $X$ there pass $|X|t$ tangents to $\st$, each of them intersects $\ell_1$ and $\ell_2$ off $\st$, hence $|X|t\leq nm$.
These two inequalities imply $q \leq t+1+nm/t$ and $|X|=q-1-t$ in the case of equality.
\Qed

\begin{proposition}
\label{dovv}
Let $\st$ be a $t$-semiarc in $\Pi_q$. If $\st$ has at least two $(q-t)$-secants and $q>2t+3$, then $\st$ is of $V_t^\circ$ type.
If $\st$ has at least two $(q-t+1)$-secants, then $\st$ is of $V_t^\bullet$ type.
\end{proposition}
\proof
If $\st$ has at least two $(q-t)$-secants with their common point in $\st$, then
Lemma \ref{j1} implies $q\leq t+1 + (t+1)^2/t = 2t+3+1/t$. If $q>2t+3$, then this is only possible when $t=1$ and $q=6$,
but there is no projective plane of order 6. Hence the common point of the $(q-t)$-secants is not contained in $\st$, which means that
$\st$ is of $V_t^\circ$ type. The proof of the second statement is straightforward.
\Qed

\begin{proposition}
\label{le2}
Let $\st$ be a $t$-semiarc of $V_t^\circ$ type in $\Pi_q$. Then the following hold.
\begin{enumerate}[(a)]
	\item $|\st| \neq 2q-2t+1$.
	\item If $t=2$, then $\st$ is a $V_2$-configuration or $|\st|=2q-2$ and $\st=(\ell_1\cup\ell_2) \triangle \Pi_2$, where $\ell_1$ and $\ell_2$ are two lines in $\Pi_2$, that is a Fano subplane contained in $\Pi_q$.
	\item If $t>1$, then $|\st| \leq 2q-t$.
\end{enumerate}
\end{proposition}
\proof
Let $\st$ be a $t$-semiarc of $V_t^\circ$ type and let $\ell_1$ and $\ell_2$ be two $(q-t)$-secants of $\st$ such that $P:=\ell_1 \cap \ell_2$ is not contained in $\st$. Denote the points of $\ell_1 \setminus (\st \cup P)$ by $A_1,\dots, A_t$, the points of $\ell_2 \setminus (\st \cup P)$ by $B_1, \dots, B_t$.
Let $X=\st \setminus(\ell_1 \cup \ell_2)$ and define the line set $\cL:=\{A_iB_j : 1 \leq i,j \leq t\}$ of size $t^2$.
Through each point $Q \in X$ there pass exactly $t$ lines of $\cL$, otherwise there would be an index $i \in \{1,2,\dots,t\}$ for
which the line $QA_i$ meets $\ell_2$ in $\st$. But then there would be at most $t-1$ tangents to $\st$ through the point $QA_i \cap \ell_2$, a contradiction. 

Suppose, contrary to our claim, that $X$ consists of a unique point denoted by $Q$. Then $Q$ would have $t+1$ tangents: the $t$ lines of $\cL$ that pass through $Q$ and the line $PQ$.

If $t=2$, then exactly two of the points of $\Pi_q \setminus (\ell_1 \cup \ell_2)$ are contained in two lines of $\cL$.
These are $Q_1:=A_1B_1 \cap A_2B_2$ and $Q_2:=A_1B_2 \cap A_2B_1$. Since $|X|>1$, we have $X=\{Q_1, Q_2\}$.
If $P$ was not collinear with $Q_1$ and $Q_2$, then $PQ_i$ would be a third tangent to $\st$ at $Q_i$, for $i=1,2$.
It follows that the point set $\Pi_2:=\{P,A_1,A_2,B_1,B_2,Q_1,Q_2\}$ is a Fano subplane in $\Pi_q$.

To prove (c), define $Y\subseteq X$ as $Y:=\{A : A\in X, |AP\cap \st|=1\}$.
The line set $\cL$ contains $|Y|(t-1)$ tangents through the points of $Y$ and $(|X|-|Y|)t$
tangents through the points of $X\setminus Y$, hence
\begin{equation}
\label{delta}
|X|(t-1) \leq |X|t-|Y| = |\cL|-\delta \leq t^2,
\end{equation}
where $\delta$ denotes the number of non-tangent lines in $\cL$.
Because of (b), we may assume $t>2$, hence $|X|\leq t^2/(t-1)< t+2$ follows.
To obtain a contradiction, suppose that $|X|=t+1$. If this is the case, then (\ref{delta}) implies $t \leq |Y|$.
If $|Y|=t$, then $X\setminus Y$ consists of a unique point, but this contradicts the definition of $Y$.
If $|Y|=t+1$, then $X=Y$ and through each point of $X$ there pass a non-tangent line, which is in $\cL$.
Thus if $\delta=1$, then the points of $X$ are contained in a line $\ell\in \cL$.
We may assume that $\ell=A_tB_t$. Then we can find $2(t-1)$ other non-tangent lines in $\cL$, these are $A_tB_i$ and $B_tA_i$ for $i=1,2,\ldots, t-1$.
On the other hand $\delta>1$ contradicts (\ref{delta}) and this contradiction proves $|X|\leq t$.
\Qed

The following result shows some kind of stability of semiarcs containing a $V_t$-configuration.

\begin{theorem}
\label{t1}
Let $\st$ be a $t$-semiarc in $\Pi_q$, $t<q$, and suppose that there exist two lines, $\ell_1$ and $\ell_2$, such that $|\ell_1 \setminus (\st\cup \ell_2)|=n$ and $|\ell_2 \setminus (\st\cup \ell_1)|=m$.
\begin{enumerate}
	\item If $\ell_1 \cap \ell_2 \notin \st$, $t>1$ and $q > \min\{n,m\}+2nm/(t-1)$, then $\st$ is of $V_t^\circ$ type.
	\item If $\ell_1 \cap \ell_2 \in \st$ and $q > \min\{n,m\}+nm/t$, then $t=(q-1)/2$, $|\st|=3(q+1)/2$ and $\st$ is of $V_t^\bullet$ type.
\end{enumerate}
We have $n=m=t$ in both cases.
\end{theorem}
\proof
In part 1, we have $n,m \geq t-1$, thus $q > \min\{n,m\}+2nm/(t-1)$ implies $n,m<q-1$.
It follows that $\ell_1$ and $\ell_2$ are not tangents to $\st$, thus $n,m \geq t$ holds.
In part 2, we have $n,m\geq t$, hence the assumption implies $n,m<q$ or, equivalently, the lines $\ell_1$ and $\ell_2$ are not tangents to $\st$.
First we show $n=m=t$ in both cases. From this, part 1 follows immediately.

We may assume $m\geq n \geq t$ and suppose, contrary to our claim, that $m\geq t+1$.
Denote by $P$ the intersection of $\ell_1$ and $\ell_2$.
Let $\cN=\{N_1,N_2,\ldots N_{q-n}\}$ be the set of points of $(\ell_1 \setminus P) \cap \st$ and $\cM=\{M_1,M_2,\ldots M_m\}$ be the set of points of $\ell_2 \setminus (\st \cup P)$. Let $X=\st \setminus (\ell_1\cup \ell_2)$.
Through each point $N_j \in \cN$ there pass exactly $m-t$ non-tangent lines that intersect $\ell_2$ in $\cM$.
Each of these lines contains at least one point of $X$. Denote the set of these points by $X(N_j)$.
Then we have the following:

\begin{itemize}
	\item $|X(N_i)|\geq m-t$, for $i=1,2,\ldots,q-n$,
	\item $X\supseteq\cup_{i=1}^{q-n} X(N_i)$,
	\item if $P\notin \st$, then each point of $X$ is contained in at most $m-t+1$ point sets of \\
	$\{X(N_1),\ldots,X(N_{q-n})\}$,
	\item if $P\in \st$, then each point of $X$ is contained in at most $m-t$ point sets of \\
	$\{X(N_1),\ldots,X(N_{q-n})\}$.
\end{itemize}
In part 1, we have the following lower bound for the size of $X$:
\begin{equation}
\label{eqqq}
\frac{(q-n)(m-t)}{m-t+1} \leq |X|.
\end{equation}
On the other hand, through each point of $X$ there pass at least $t-1$ tangents that intersect both $\ell_1 \setminus (\st \cup P)$ and $\cM$.
Hence we have:
\begin{equation}
\label{eqqqq}
|X|\leq \frac{nm}{t-1}.
\end{equation}
Summarizing these two inequalities we get:
$$q \leq n + \frac{nm}{t-1}+\frac{nm}{(m-t)(t-1)}\leq n + \frac{2nm}{t-1},$$ that is a contradiction.

In part 2, observe that Lemma \ref{j1} and $q>\min\{m,n\}+nm/t$ together imply $n=t$.
If $m\geq t+1$, then similarly to (\ref{eqqq}) and (\ref{eqqqq}), we get $(q-t)(m-t)/(m-t) \leq |X|$ and $|X|\leq mt/t$ respectively.
These two inequalities imply $q \leq t+m$, contradicting our assumption $q>\min\{n,m\}+nm/t=t+m$, hence $m=t$ follows.
If $n=m=t$, then Lemma \ref{j1} implies $q\leq 2t+1$ while our assumption yields $q>2t$, thus $t=(q-1)/2$.
Since in this case there is equality in Lemma \ref{j1}, we have $|\st|=3q-3t=3(q+1)/2$.
\Qed

Let $\alpha_{n,m}$ and $\beta_{n,m}$ denote the lower bounds on $q$ in part 1 and in part 2 of Theorem \ref{t1}, respectively.
The following example shows that the weaker assumptions $\alpha_{n,m}<3q$ and $\beta_{n,m}<2q$, respectively, do not imply
the existence of a $V_t$-configuration contained in the semiarc.
%$n=m=t$.

\begin{example}
We give two examples for $t$-semiarcs, $\st$, such that they do not contain a $V_t$-configuration and there exist two lines, $\ell_1$ and $\ell_2$, with $\ell_1\setminus (\ell_2 \cup \st)=t$ and $\ell_1\setminus (\ell_2 \cup \st)=t+1$.
To do this, choose a conic $\cC$ in $\Pi_s$, that is a projective plane of order $s>3$.
Let $Q_1$ and $Q_2$ be two points of $\cC$ and proceed as follows.
\begin{enumerate}
	\item Let $\ell_i$ be the tangent of $\cC$ at the point $Q_i$, for $i=1,2$, and denote $\ell_1 \cap \ell_2$ by $P$.
Take a point $Z \in Q_1Q_2$ such that $PZ$ is a secant of $\cC$. Then
$\so_0:=(\ell_1 \cup \ell_2 \cup \cC \cup \{Z\})\setminus \{P,Q_2\}$ is a point set without tangents.
Now, if $\Pi_s$ is contained in $\Pi_q$, then $\so_0\subset \Pi_s$ is a $t$-semiarc in $\Pi_q$, with $t=q-s$.
We have $\ell_1 \cap \ell_2 \notin \st$ and
$$\alpha_{t,t+1}=(q-s) + 2\frac{(q-s+1)(q-s)}{q-s-1}<3q.$$
	\item Let $\ell_1$ be the tangent of $\cC$ at $Q_1$ and let $\ell_2$ be the line $Q_1Q_2$.
Take a point $Z\in \ell\setminus(\ell_1\cup\ell_2)$, where $\ell$ denotes the tangent of $\cC$ at $Q_2$. Then
$\so_0:=(\ell_1 \cup \ell_2 \cup \cC \cup \{Z\})\setminus \{Q_2\}$ is a point set without tangents.
As before, if $\Pi_s$ is contained in $\Pi_q$, then $\so_0\subset \Pi_s$ is a $t$-semiarc in $\Pi_q$, with $t=q-s$.
We have $\ell_1 \cap \ell_2 \in \st$ and
$$\beta_{t,t+1}=(q-s)+\frac{(q-s+1)(q-s)}{q-s}<2q.$$
\end{enumerate}
\end{example}

%If $q>4$, then a vertexless triangle with one point deleted from one of its sides is a semioval $\so_1$ in $\Pi_q$, which has a $(q-1)$-secant
%and a $(q-2)$-secant intersecting each other not in $\so_1$. This example shows that the assumption $t>1$ in part 1 of Theorem \ref{t1} is necessary to
%conclude $n=m=t$. However, this example is of $V_1^\circ$ type.

The next example is due to Suetake and it shows that when $t=1$, then there is no analogous result for part 1 of Theorem \ref{t1}.

\begin{example}[\cite{suetake}, Example 3.3]
Let $A$ be a proper, not empty subset of $\GF(q)\setminus \{0\}$, such that $A=-A:=\{-a : a\in A\}$ and $|A|\geq 2$.
Let $B=\GF(q)\setminus (A\cup \{0\})$ and define the following set of points in $\PG(2,q)$:
$$\so_1:=\{(0,a,1),(b,0,1),(c,c,1),(m,1,0) : a\in A, b\in B, c\in \GF(q)\setminus\{0\}, m\in \GF(q)\setminus\{0,1\}\}.$$
Then $\so_1$ is a semioval with a $(q-1)$-secant, $X=Y$, and a $(q-2)$-secant, $Z=0$, intersecting each other not in $\so_1$.
Also, $\so_1$ is not of $V_1^\circ$ type.
\end{example}

When $A=\GF(q)\setminus \{0\}$ in the above example, then $\so_1$ is a vertexless triangle with one point deleted from one of its sides.
This example exists also in non-Desarguesian planes, but it is a semioval of $V_1^\circ$ type.

%If $q>4$, then a vertexless triangle with one point deleted from one of its sides is a semioval $\so_1$ in $\Pi_q$, which has a $(q-1)$-secant
%and a $(q-2)$-secant intersecting each other not in $\so_1$. This example exist also in non-Desarguesian planes, however, it is of $V_1^\circ$ type.

Semiarcs that properly contain a $V_t$-configuration exist in $\Pi_q$ whenever $\Pi_q$ contains a subplane.
Some of the following examples were motivated by an example due to Korchm\'aros and Mazzocca (see $\cite{nuclei}$, pg. 64).

\begin{example}
\label{KM}
Let $\Pi^0, \Pi^1, \dots, \Pi^{s-1}$ be subplanes of $\Pi^s:=\Pi_q$ such that $\Pi^{i-1} \subset \Pi^i$ for $i=1,\ldots,s$.
Denote by $r$ the order of $\Pi^0$ and let $\ell_1$ and $\ell_2$ be two lines in this plane. Let $P=\ell_1\cap \ell_2$ and set
$$S(0):=(\ell_1\cup \ell_2)\cap (\Pi^0 \setminus P),~S(j):=(\ell_1\cup \ell_2)\cap(\Pi^j \setminus \Pi^{j-1}),~\text{for}~j=1,\ldots,s.$$
By $I$ we denote a subset of $\{1,2,\ldots,s\}$. We give four examples.
\begin{enumerate}
	\item Let $\ell$ be a line in $\Pi^0$ passing through $P$ and let $Z$ be a subset of $(\ell \cap \Pi^0) \setminus \{P\}$ of size at least two.
	If $I$ is not empty, then $\st:=\cup_{j\in I} S(j)\cup Z$ is a $t$-semiarc of $V_t^\circ$ type with $t=q-\frac{1}{2}\sum_{j\in I}|S(j)|$.
	\item Let $\ell$ be a line in $\Pi^0$ that does not pass through $P$ and let $Z$ be a subset of $(\ell \cap \Pi^0) \setminus (\ell_1 \cup \ell_2)$ of size
	at least two. If $I$ is not empty, then $\st:=\cup_{j\in I} S(j)\cup Z$ is a $t$-semiarc of $V_t^\circ$ type with $t=q-\frac{1}{2}\sum_{j\in I}|S(j)|$.
	\item Let $Z$ be a subset of $\Pi^0 \setminus (\ell_1\cup \ell_2)$ such that there is no line in $\Pi^0$ passing through $P$
	and meeting $Z$ in exactly one point.
	If $I$ is a proper subset of $\{1,2,\ldots,s\}$, then $\st:=\cup_{j\in I} S(j)\cup Z \cup S(0)$ is a $t$-semiarc of $V_t^\circ$ type with $t=q-r-\frac{1}{2}\sum_{j\in I}|S(j)|$.
	\item Let $Z$ be a subset of $\Pi^0 \setminus (\ell_1\cup \ell_2)$ such that for each line $\ell\neq \ell_1,\ell_2$ through $P$, $\ell$ is a line in $\Pi^0$, we have	$|\ell \cap Z|\geq 1$. Then $\st:=\{P\}\cup S(0)\cup Z$ is a $t$-semiarc of $V_t^\bullet$ type with $t=q-r$.
\end{enumerate}
\end{example}

\section{Semiarcs containing a $V_t$-configuration in $\PG(2,q)$}
\label{sec:1.3}
\indent

In this section our aim is to characterize $t$-semiarcs containing a $V_t$-configuration in $\PG(2,q)$.
We will need the following definition.

\begin{definition}
Let $\ell_1$ and $\ell_2$ be two lines in a projective plane and let $P$ denote their common point.
We say that $X_1\subseteq \ell_1\setminus P$ and $X_2 \subseteq \ell_2 \setminus P$ are two perspective point sets if
there is a point $Q$ such that each line through $Q$ intersects both $X_1$ and $X_2$ or intersects none of them.
In other words, there is a perspectivity which maps $X_1$ onto $X_2$.
\end{definition}

\begin{lemma}
\label{le}
Let $\st$ be a $t$-semiarc in $\Pi_q$ and suppose that $(\ell_1 \triangle \ell_2)\cap \st$ is a $V_t$-configuration for some lines $\ell_1$ and $\ell_2$.
If $\st \nsubseteq \ell_1 \cup \ell_2$, then $\st \cap (\ell_1\setminus \ell_2)$ and $\st\cap (\ell_2 \setminus \ell_1)$ are perspective point sets and each point of $\st\setminus(\ell_1 \cup \ell_2)$ is the centre of a perspectivity which maps $\st \cap (\ell_1\setminus \ell_2)$ onto $\st\cap (\ell_2 \setminus \ell_1)$.
\end{lemma}
\proof
Let $X=\st \setminus (\ell_1 \cup \ell_2)$ and $X_i=\st \cap (\ell_i\setminus \ell_j)$, for $\{i,j\}=\{1,2\}$.
For each $Q\in X$, if there was a line $\ell$ through $Q$ intersecting $X_i$ but $X_j$, 
then the point $\ell \cap X_i \in \st$ would have at most $t-1$ tangents.
This shows that each point of $X$ is the centre of a perspectivity which maps $X_1$ onto $X_2$.
If $\st \nsubseteq \ell_1 \cup \ell_2$, then $X$ is not empty, hence $X_1$ and $X_2$ are perspective point sets.
\Qed

The following theorem characterizes perspective point sets in $\PG(2,q)$.
This result was first published by Korchm\'aros and Mazzoca in \cite{nuclei} but we will
use the notation of \cite{napoli} by Bruen, Mazzocca and Polverino.

\begin{theorem}[\cite{napoli}, Result 2.2, Result 2.3, Result 2.4, see also \cite{nuclei}]
\label{persp}
Let $\ell_1$ and $\ell_2$ be two lines in $\PG(2,q)$, $q=p^r$, and let $P$ denote their common point.
Let $X_1\subseteq \ell_1\setminus P$ and $X_2 \subseteq \ell_2 \setminus P$ be two perspective point sets.
Denote by $U$ the set of all points which are centres of a perspectivity mapping $X_1$ onto $X_2$.
Using a suitable projective frame in $\PG(2,q)$, there exist an additive subgroup $B$ of $\GF(q)$ and a multiplicative subgroup
$A$ of $\GF(q)$ such that:
\begin{enumerate}[(a)]
	\item $B$ is a subspace of $\GF(q)$ of dimension $h_1$ considered as a vectorspace over a subfield $\GF(q_1)$ of $\GF(q)$ with $q_1=p^d$
	and $d|r$. This implies that $B$ is an additive subgroup of $\GF(q)$ of order $p^h$ with $h=dh_1$.
	\item $A$ is a multiplicative subgroup of $\GF(q_1)$ of order $n$, where $n|(p^d-1)$. In this way, $B$ is invariant under $A$, i.e.
	$B=AB := \{ab : a\in A,b\in B\}$.
	\item If $G_i$ denotes the full group of affinities of $\ell_i\setminus P$	preserving the set $X_i$, $i=1,2$, then
	$G_1 \cong G_2 \cong G=G(A,B)=\{g : g(y)=ay+b,a\in A,b\in B\}\leq \Sigma$, where $\Sigma$ is the full affine group on the line $\AG(1,q)$.
	\item $X_i$ is a union of orbits of $G_i$ on $\ell_i\setminus P$, $i=1,2$, and $|U|=|G|=np^h$.
	\item For every two integers $n, h$, such that $n|(p^d-1)$ and $d|gcd(r,h)$, there exists in $\Sigma$
	a subgroup of type $G=G(A,B)$ of order $np^h$, where $A$ and $B$ are multiplicative and	additive subgroups of $\GF(q)$ of order $n$ and $p^h$, respectively.
	\item $G$ has one orbit of length $p^h$ on $\AG(1,q)$, namely $B$, and $G$ acts regularly on the remaining orbits, say $O_1,O_2,\ldots, O_m$,
	where	$$m=\frac{q-p^h}{np^h}=\frac{p^{r-h}-1}{n}.$$
\end{enumerate}
In the sequel we denote by $B^i$ the orbit of $G_i$ on $\ell_i \setminus P$ corresponding to $B$ and by $O_1^i,O_2^i,\ldots,O_m^i$ the remaining orbits, for $i=1,2$.
With this notation $B^1$ is the image of $B^2$ under the perspectivities with centre in $U$ and also $O_j^2$ is the image of $O_j^1$ for $j=1,2,\ldots, m$ and vice versa.
\begin{enumerate}[(a)]
  \setcounter{enumi}{6}
	\item $B^1\subseteq X_1$ if and only if $B^2 \subseteq X_2$ and the same holds for the other orbits, i.e. $O_j^1 \subseteq X_1$ if and only if
	$O_j^2 \subseteq X_2$, for $j=1,2,\ldots, m$.
	\item If a line $\ell$ not through $P$ meets $U$ in at least two points, then $\ell$ intersects both $B^1$ and $B^2$.
\end{enumerate}
Exactly one of the following cases must occur.
\begin{enumerate}
	\item Both $A$ and $B$ are trivial. Then $U$ consists of a singleton.
	\item $A$ is trivial and $B$ is not trivial. Then $U$ is a set of $p^h$ points all collinear with the point $P$.
	\item $B$ is trivial and $A$ is not trivial. Then $U$ is a set of $n$ points on a line not through $P$.
	\item $A$ and $B$ are the multiplicative and the additive group, respectively, of a subfield $\GF(p^h)$ of $\GF(q)$. Then
	$$U\cup B^1\cup B^2\cup \{P\}= \PG(2,p^h).$$
	\item None of the previous cases occur. Then $U$ is a point set of size $np^h$ and of type $(0, 1, n, p^h)$, i.e. 0, 1, $n$, $p^h$ are the only intersection numbers of $U$ with respect to the lines in $\PG(2,q)$. In addition, using the fact that $|U|=np^h$,
	\begin{itemize}
	\item there are exactly n lines intersecting U in exactly $p^h$ points and they are all concurrent at the common point $P$ of
	$\ell_1$ and $\ell_2$,
	\item each line intersecting U in exactly n points meets both $B^1$ and $B^2$.
\end{itemize}
\end{enumerate}
\end{theorem}

%Semiarcs contained in the union of three lines such that their common point is contained in the semiarc were described is $\cite{csk}$.
%Since $(q-2)$-semiarcs satisfy this property, we have the following result:
%The following result followss from [14], Proposition 3.1.

\begin{lemma}[\cite{csk}, Proposition 3.1]
\label{notes}
If $\st$ is a $(q-2)$-semiarc in $\Pi_q$, then it is one of the following three configurations:
four points in general position, the six vertices of a complete quadrilateral, or a Fano subplane.
\end{lemma}

In the next theorems we will use the notation of Theorem $\ref{persp}$.

\begin{theorem}
\label{thm}
Let $\st$ be a $t$-semiarc in $\PG(2,q)$, $q=p^r$, and suppose that $(\ell_1 \triangle \ell_2)\cap \st$ is a $V_t$-configuration for some lines $\ell_1$ and $\ell_2$. To avoid trivial cases, suppose that $\st\nsubseteq \ell_1 \cup \ell_2$.

Let $X_i=\ell_i\cap \st$, for $i=1,2$, and let $X = \st \setminus (\ell_1 \cup \ell_2)$. Also let $P=\ell_1 \cap \ell_2$.
Because of Lemma \ref{le} we have that $X_1$ and $X_2$ are perspective point sets and $X\subseteq U$,
where $U$ is the set of all points which are centres of a perspectivity mapping $X_1$ onto $X_2$.
Choose a suitable coordinate system as in Theorem \ref{persp} and suppose that the size of $G=G(A,B)$ is $np^h$, i.e.
$|A|=n$ and $|B|=p^h$, where $A$ and $B$ are the multiplicative and the additive subgroup of $\GF(q)$ associated to the perspective point sets $X_1$ and $X_2$.

\begin{enumerate}[(I)]
\item If $P\notin \st$, i.e. $\st$ is of $V_t^\circ$ type, then one of the following holds.
\begin{enumerate}[(i)]
	\item $X$ is contained in a line through $P$ that meets $U$ in $p^h$ points, $h\geq 1$, and we have $2 \leq |X| \leq p^h$,
	\item $X$ is contained in a line not through $P$ that meets $U$ in $n\geq 2$ points and we have $2 \leq |X| \leq n$,
	\item $|X|\geq 2$ and $X$ is a subset of $U$ such that there is no line through $P$ that meets $X$ in exactly one point.
\end{enumerate}
\end{enumerate}

In the first two cases $X_i=\cup_{j\in I}O^i_j$ for some not empty subset $I \subseteq \{1,2,\ldots,m\}$ and for $i=1,2$.
We have $t=q-knp^h$, where $k=|I|$ and $1 \leq k \leq m$, where $m=(p^{r-h}-1)/n$.

In the third case $X_i=\cup_{j\in I}O^i_j \cup B^i$ for some proper subset $I \subset \{1,2,\ldots,m\}$ and for $i=1,2$.
We have $t=q-knp^h-p^h$, where $k=|I|$, $h \geq 1$ and $0 \leq k \leq m-1$.

\begin{enumerate}[(I)]
\setcounter{enumi}{1}
\item If $P\in \st$, i.e. $\st$ is of $V_t^\bullet$ type, then one of the following holds.
\begin{enumerate}[(i)]
	\item $\st$ consists of the six vertices of a complete quadrilateral or $\st$ is a Fano subplane.
	We have $t=q-2$ in both cases.
%	\item $\st$ is a Fano subplane, in this case $t=q-2$. 
	\item $\ell_1$ and $\ell_2$ are lines in the subplane $\PG(2,p^h)$ and
	$$\st=\PG(2,p^h) \cap (\ell_1 \cup \ell_2)\cup X,$$
  where $X$ is a subset of $\PG(2,p^h) \setminus (\ell_1\cup \ell_2)$ such that for each line $\ell\neq \ell_1,\ell_2$ through $P$, $\ell$ is a line in $\PG(2,p^h)$, we have	$|\ell \cap X|\geq 1$. In this case $t=q-p^h$.
	\item $\st$ is projectively equivalent to the following set of $3(n+1)$ points:
	$$\st:=\{(a,0,1),(0,-a,1),(a,1,0): a \in A\}\cup\{(1,0,0),(0,1,0),(0,0,1)\},$$
%	where $A$ is a multiplicative subgroup of $\GF(q)$ of order $n$.
	In this case $t=q-1-n$, where $n\mid q-1$.
\end{enumerate}
\end{enumerate}
The converse is also true, if $X_1$ and $X_2$ are perspective point sets and $X$ is as in one of the three cases in (I), then $X\cup X_1 \cup X_2$
is a $t$-semiarc of $V_t^\circ$ type. If $\st$ is as in one of the three cases in (II), then $\st$ is a $t$-semiarc of $V_t^\bullet$ type.
\end{theorem}

\proof
We begin by proving (I).
First assume $B^1 \subseteq \ell_1 \setminus X_1$. Then Theorem \ref{persp} (g) implies  $B^2 \subseteq \ell_2 \setminus X_2$.
Suppose that there exist three non-collinear points in $X$, say $L, M$ and $N$.
Then between the lines $LM, LN$ and $MN$ there are at least two, say $LM$ and $LN$, not through $P$.
Theorem \ref{persp} (h) and $X\subseteq U$ imply that these two lines intersect both $B^1$ and $B^2$.
But then through $L$ there pass at most $t-1$ tangents, a contradiction.
It follows that $X$ is contained in a line and hence it is as in one of our first two cases.
The condition $|X|\geq 2$ comes from Proposition \ref{le2} (a).

Now asume $B^1 \subseteq X_1$ and hence $B^2 \subseteq X_2$. In this case for every two points $M,N\in X$, the line $MN$ intersects
$\ell_i$ in $X_i$, for $i=1,2$. Thus the number of tangents through a point $L\in X$ is $t$ if and only if the line $LP$ contains at least
one other point of $X$. Case 3 of Theorem \ref{persp} shows that this is not possible when $B$ is trivial, i.e. when $h=0$.
Hence $X$ is as in our third case.

Now we prove (II).
First assume $B^1 \subseteq \ell_1 \setminus X_1$ and hence $B^2 \subseteq \ell_2 \setminus X_2$.
Suppose that there exist two points in $X$, say $M$ and $N$, not collinear with $P$.
Then the line $MN$ intersects $\ell_1$ and $\ell_2$ not in $\st$.
But then the number of tangents through $M$ is at most $t-1$, a contradiction.
Thus $X$ is contained in a line through $P$ and through $P$ there pass exactly $q-2$ tangents.
So $\st$ is a $(q-2)$-semiarc. According to Lemma \ref{notes}, $\st$ is as in (II)(i).

Now assume $B^1 \subseteq X_1$ and hence $B^2 \subseteq X_2$.
In this case $t=q-knp^h-p^h$ for some $k\in \{0,1,\ldots,m-1\}$, where $m$ is the number of orbits of $G$ of size $np^h$ on $\AG(1,q)\setminus B$.
Since $P$ has exactly $t$ tangents, there are $q+1-t$ non-tangent lines through $P$. According to Theorem \ref{persp}, we have $q-1-t \leq n$ and
hence $knp^h+p^h-1 \leq n$. We distinguish two subcases.

If $h>0$, then $n|p^h-1$ implies $n \leq p^h-1$ and hence $knp^h=0$. This occurs only if $k=0$ and $n=p^h-1$.
But $n$ divides also $p^d-1$, where $d|h$ and $B$ is a subspace over the field $\GF(p^d)$. This implies $d=h$, thus $B$ is a subfield and $U$ is
as in case 4 of Theorem \ref{persp}. This is only possible if $\st$ is as in our second case.

If $h=0$, then $kn \leq n$ and $U$ is as in case 3 of Theorem \ref{persp}. If $k=0$, then $t=q-1$, which we excluded.
Thus we have $k = 1$ and $t=q-n-1$. This occurs only if $\st$ is as in our third case (see \cite{napoli}, pg. 56--57).
\Qed

\begin{theorem}
Let $\st$ be a $t$-semiarc of $V_t^\circ$ type in $\PG(2,q)$, $q=p^r$. Then the following hold.
	\begin{enumerate}[(a)]
		\item If $\gcd(q,t)=1$ and $\gcd(q-1,t-1)=1$, then $\st$ is a $V_t$-configuration.
		\item If $\gcd(q,t)=1$, then $\st$ is contained in a vertexless	triangle.
		\item If $\gcd(q-1,t)=1$, then $\st$ is contained in a vertexless	triangle or in the union of three concurrent lines without their common point.
	\end{enumerate}
\end{theorem}

\proof
We have $p^h|t$ in all three cases of Theorem \ref{thm} (I), where $p^h$ is the size of $B$.
Hence $\gcd(q,t)=1$ implies $p^h=1$, i.e. $h=0$. This occurs only in the second case of Theorem \ref{thm} (I) and this proves (b).

In the first two cases of Theorem \ref{thm} (I) we have $n|(t-1)$ and hence also $n|\gcd(q-1,t-1)$, where $n$ is the size of $A$.
We have seen previously that $\gcd(q,t)=1$ can hold only in the second case of Theorem \ref{thm} (I).
But in that case we have $n\geq 2$, which is a contradiction when $\gcd(q-1,t-1)=1$. This proves (a).

If $\st$ is as in one of the first two cases of Theorem \ref{thm} (I), then we are done.
So to prove (c), it is enough to consider Theorem \ref{thm} (I)(iii).
In this case $t=(q-1)-nkp^h-(p^h-1)$ and hence $n|\gcd(q-1,t)$.
If $\gcd(q-1,t)=1$, then $n=1$, i.e. $A$ is trivial.
If this happens, then case 2 of Theorem \ref{persp} implies that $\st$ is contained in the union of three concurrent lines without their common point.
\Qed

{\bf Acknowledgement.}
The author is very grateful for the advices of Prof. G\'{a}bor Korchm\'{a}ros and Prof. Tam\'as Sz\H{o}nyi.

\begin{flushleft}
Bence Csajb\'{o}k \\
Department of Mathematics, Informatics and Economics\\
University of Basilicata \\
Campus Macchia Romana, via dell'Ateneo Lucano \\
I-85100 Potenza, Italy \\
e-mail: {\sf bence.csajbok@unibas.it}
\end{flushleft}

\end{document}